\documentclass{article}

\usepackage[letterpaper]{geometry}
\usepackage{hyperref}
\usepackage{titling}

\usepackage[square, numbers]{natbib} 
\usepackage{amsmath}
\usepackage{amsthm}
\usepackage{bm}
\usepackage{amssymb}
\usepackage[mathscr]{eucal}
\usepackage{graphicx}
\usepackage{enumerate}
\usepackage[linesnumbered,boxed,titlenumbered,ruled,nofillcomment,noend]{algorithm2e}
\usepackage{booktabs}
\usepackage{scalerel}
\usepackage[dvipsnames]{xcolor}

\hypersetup{
	colorlinks=true,
	linkcolor=Blue,
	urlcolor=Blue,
	citecolor=Blue,
}

\theoremstyle{plain}
\newtheorem{theorem}{Theorem} \numberwithin{theorem}{section}

\theoremstyle{definition}

\newtheorem{definition}[theorem]{Definition}

\newcommand{\Ae}{\boldsymbol{\mathscr{A}}}

\newcommand{\Xe}{\boldsymbol{\mathscr{X}}}
\newcommand{\Ye}{\boldsymbol{\mathscr{Y}}}

\newcommand{\F}{\textup{F}}

\newcommand{\Phibf}[0]{\bm{\Phi}}

\newcommand{\Abf}[0]{\bm{A}}
\newcommand{\Bbf}[0]{\bm{B}}
\newcommand{\Cbf}[0]{\bm{C}}

\newcommand{\Gbf}[0]{\bm{G}}

\newcommand{\Mbf}[0]{\bm{M}}

\newcommand{\Sbf}[0]{\bm{S}}

\newcommand{\Xbf}[0]{\bm{X}}
\newcommand{\Ybf}[0]{\bm{Y}}

\newcommand{\ebf}[0]{\bm{e}}

\newcommand{\jbf}[0]{\bm{j}}

\newcommand{\pbf}[0]{\bm{p}}

\newcommand{\xbf}[0]{\bm{x}}
\newcommand{\ybf}[0]{\bm{y}}

\newcommand\Pb{\mathbb{P}}

\newcommand\Rb{\mathbb{R}}

\DeclareMathOperator*{\argmin}{arg\,min}

\DeclareMathOperator*{\Ind}{\operatorname{Ind}}
\newcommand{\defeq}{\stackrel{\text{\tiny \textnormal{def}}}{=}}

\newcommand{\rank}[0]{\operatorname{rank}}
\newcommand{\noiter}{\textup{\#it}}

\newcommand*{\OPT}{\textup{OPT}}

\newcommand{\trace}[0]{\operatorname{trace}}

\newcommand{\TN}[0]{\operatorname{TN}}

\DeclareMathOperator*{\startimes}{\scalerel*{\circledast}{\sum}}

\newtheorem{Example}{Example}[section]

\begin{document}

\title{Sampling-Based Decomposition Algorithms for Arbitrary Tensor Networks}

\author{
	Osman Asif Malik\thanks{Lawrence Berkeley National Laboratory, \href{mailto:oamalik@lbl.gov}{\texttt{oamalik@lbl.gov}}}
	\and Vivek Bharadwaj\thanks{UC Berkeley, \href{mailto:vivek_bharadwaj@berkeley.edu}{\texttt{vivek\_bharadwaj@berkeley.edu}}}
	\and Riley Murray\thanks{UC Berkeley, \href{mailto:rjmurray@berkeley.edu}{\texttt{rjmurray@berkeley.edu}}}}

\predate
\postdate
\date{}

\begin{NoHyper}
\maketitle
\end{NoHyper}

\begin{abstract}
	We show how to develop sampling-based alternating least squares (ALS) algorithms for decomposition of tensors into \emph{any} tensor network (TN) format. 
	Provided the TN format satisfies certain mild assumptions, resulting algorithms will have \emph{input sublinear} per-iteration cost.
	Unlike most previous works on sampling-based ALS methods for tensor decomposition, the sampling in our framework is done according to the \emph{exact} leverage score distribution of the design matrices in the ALS subproblems.
	We implement and test two tensor decomposition algorithms that use our sampling framework in a feature extraction experiment where we compare them against a number of other decomposition algorithms.
\end{abstract}

\section{Introduction} \label{sec:intro}

Tensor decomposition has emerged as an important tool in data mining and machine learning \citep{papalexakis2016TensorsData, cichocki2016TensorNetworksPart1, cichocki2017TensorNetworksPart2, ji2019SurveyTensor}.
Applications in data mining include network analysis, web mining, topic modeling and recommendation systems.
Tensor decomposition is used widely in machine learning for things like parameter reduction in neural networks, understanding of deep neural network expressiveness, supervised learning and feature extraction. 

Due to the multidimensional nature of tensors, they are inherently plagued by the curse of dimensionality.
For example, representing a tensor $\Xe\in \Rb^{I \times \cdots \times I}$ with $N$ modes requires $I^N$ numbers.
This exponential dependence on $N$ makes its way into algorithms for computing tensor decompositions.
Alternating least squares (ALS) is arguably the most popular and successful approach for computing a wide range of tensor decompositions.
When decomposing a tensor $\Xe$, each iteration of ALS involves solving a sequence of least squares problems for which the entries of $\Xe$ feature as the dependent variables.
The per-iteration cost of ALS therefore naturally inherits the exponential dependence on $N$.

A large number of papers have sought to reduce the cost of tensor decomposition by leveraging techniques from randomized numerical linear algebra.
One particularly interesting line of work \citep{cheng2016SPALSFast, larsen2022PracticalLeverageBased, fahrbach2022SubquadraticKronecker, malik2021SamplingBasedMethod, malik2022MoreEfficient} seeks to construct ALS algorithms with a per-iteration cost which is \emph{sublinear} in the number of tensor entries, i.e., $o(I^N)$.
Since any algorithm considering all entries of $\Xe$ immediately incurs a cost of $\Omega(I^N)$, these works have all resorted to \emph{sampling-based} techniques.
More precisely, they all sample the ALS subproblems according to the leverage score distribution (or an approximation thereof).
This is done efficiently by taking advantage of the special structure of the design matrices in the ALS subproblems for tensor decomposition.

The previous works discussed above develop methods for specific tensor decompositions: 
The CP decomposition in \citep{cheng2016SPALSFast, larsen2022PracticalLeverageBased, malik2022MoreEfficient}, Tucker decomposition in \citep{fahrbach2022SubquadraticKronecker}, and the tensor ring decomposition in \citep{malik2021SamplingBasedMethod, malik2022MoreEfficient}.
In this paper, we consider \emph{all} decompositions that can be expressed in tensor network (TN) format.
The TN format allows for a very wide range of decompositions, including the CP, Tucker, tensor train and tensor ring decompositions.
The following summarizes our contributions:
\begin{itemize}
	\item We first show how to efficiently sample rows of \emph{any} tall-and-skinny matrix which is in TN format according to the \emph{exact} leverage score distribution.
	
	\item We then show how this sampling technique can be used to yield ALS algorithms with a per-iteration cost which is input sublinear for \emph{all} TN decompositions that satisfy certain mild assumptions.
\end{itemize}

The decomposition framework we present builds on the work in \citep{malik2022MoreEfficient}.
That paper is notable since it provided the first sampling-based ALS methods for CP and tensor ring decomposition with a per-iteration cost depending \emph{polynomially} on $N$; 
the dependence in earlier works was exponential \citep{cheng2016SPALSFast, larsen2022PracticalLeverageBased, malik2021SamplingBasedMethod}. 
We are able to substantially simplify the scheme in \citep{malik2022MoreEfficient} by entirely avoiding the complicated recursive sketching procedure that it relies on.
This makes implementation easier and is also what ensures that the sampling is done according to the exact leverage score distribution rather than an approximation of it.
The simplification is also what paves the way for generalization to arbitrary TN formats.

\section{Related Work} \label{sec:related-work}

\citet{cheng2016SPALSFast} develop the first ALS method for CP decomposition with an input sublinear per-iteration cost. 
Their method uses a mixture of leverage score and row-norm sampling applied to the matricized-tensor-times-Khatri--Rao product (MTTKRP) which arises as a key computational kernel in CP decomposition.
\citet{larsen2022PracticalLeverageBased} use leverage score sampling to reduce the size of each ALS subproblem.
Their method has improved theoretical guarantees compared to those in \citep{cheng2016SPALSFast} as well as improved scalability due to various practical improvements.
\citet{malik2021SamplingBasedMethod} propose a sampling-based ALS method for tensor ring decomposition which also uses leverage score sampling to reduce the size of the ALS subproblems.

All three papers \citep{cheng2016SPALSFast, larsen2022PracticalLeverageBased, malik2021SamplingBasedMethod} require a number of samples which scales exponentially with the number of input tensor modes, $N$, for performance guarantees to hold.
This translates into a per-iteration cost which is $\Omega(R^{N+1})$ for the CP decomposition \citep{cheng2016SPALSFast, larsen2022PracticalLeverageBased} and $\Omega(R^{2N+2})$ for the tensor ring decomposition \citep{malik2021SamplingBasedMethod}, where $R$ is the relevant notion of rank.
For low-rank decompositions ($R \ll I$), the cost will be input sublinear despite this exponential cost, but for higher rank it might not be (recall that, unlike in matrix decomposition, we can have $R > I$ in tensor decomposition).
\citet{malik2022MoreEfficient} develop methods for both CP and tensor ring decomposition which avoid this exponential dependence on $N$, instead improving it to a polynomial dependence.
This is achieved by sampling from a distribution much closer to the exact leverage score distribution than the previous works \citep{cheng2016SPALSFast, larsen2022PracticalLeverageBased, malik2021SamplingBasedMethod} do.
Since our work builds on and improves the scheme in \citep{malik2022MoreEfficient}, it also has a polynomial dependence on $N$ when used for CP and tensor ring decomposition (see Tables~\ref{tab:cp-complexity-comparison} and \ref{tab:tr-complexity-comparison}).

\citet{fahrbach2022SubquadraticKronecker} develop a method for efficient sampling of ridge regression problems involving Kronecker product design matrices.
They use these techniques to achieve an efficient sampling-based ALS method for regularized Tucker decomposition.

Efficient sketching of structured matrices is an active research area with applications beyond tensor decomposition.
The recent work by \citet{ma2022CostefficientGaussian} is particularly interesting since it proposes structured sketching operators for general TN matrices.
These operators take the form of TNs made up of Gaussian tensors and are therefore dense.
When used in ALS for tensor decomposition, their sketching operators therefore need to access all entries of the data tensor in each least squares solve which leads to a per-iteration cost which is at least linear in the input tensor size (and therefore exponential in $N$).

Due to space constraints, we have focused on previous works that develop sampling-based ALS methods here.
There is a large number of works that develop randomized tensor decomposition methods using other techniques.
For a comprehensive list of such works, see \citep[Sec.~2]{malik2022MoreEfficient}.
For an overview of work on structured sketching, see \citep[Sec.~1.2]{malik2022FastAlgorithms}.

\section{Preliminaries} \label{sec:preliminaries}

\subsection{General Notation}

By \emph{tensor}, we mean a multidimensional array containing real numbers.
We refer to a tensor with $N$ indices as an \emph{$N$-way tensor} and we say that it is of \emph{order} $N$ or that it has $N$ \emph{modes}.
We use bold uppercase Euler script letters (e.g., $\Xe$) to denote tensors of order 3 or greater, bold uppercase letters (e.g., $\Xbf$) to denote matrices, bold lowercase letters (e.g., $\xbf$) to denote vectors, and regular lowercase letters to denote scalars (e.g., $x$).
We denote entries of an object either in parentheses or with subscripts.
For example, $\Xe(i,j,k) = \Xe_{ijk}$ is the entry on position $(i,j,k)$ in the 3-way tensor $\Xe$, and $\xbf(i)=\xbf_i$ is the $i$th entry in the vector $\xbf$.
We use a colon to denote all entries along a particular mode.
For example, $\Xbf(i,:)$ denotes the $i$th row of the matrix $\Xbf$.
Superscripts in parentheses will be used to denote a sequence of objects.
For example, $\Ae^{(1)},\ldots,\Ae^{(M)}$ is a sequence of tensors and $\ebf^{(i)}$ is the $i$th canonical basis vector.
The values that a tensor's indices can take are referred to as \emph{dimensions}.
For example, if $\Xe = (\Xe_{ijk})_{i=1, j=1, k=1}^{I, J, K}$ then the dimensions of modes 1, 2 and 3 are $I$, $J$ and $K$ respectively.
Note that dimension \emph{does not} refer to the number of modes, which is 3 for $\Xe$.
For a positive integer $I$, we use the notation $[I] \defeq \{1, \ldots, I\}$.
For indices $i_1 \in [I_1], \ldots, i_N \in [I_N]$, the notation $\overline{i_1 \cdots i_N} \defeq 1 + \sum_{n=1}^N(i_n-1)\prod_{j=1}^{n-1} I_j$ will be used to denote the linear index corresponding to the multi-index $(i_1,\ldots,i_N)$.
The pseudoinverse of a matrix $\Abf$ is denoted by $\Abf^+$.
We use $\otimes$, $\odot$ and $\circledast$ to denote the Kronecker, Khatri--Rao and Hadamard matrix products which are defined in Section~\ref{sec:additional-notation}.

\subsection{Tensor Networks} \label{sec:tensor-networks}

\begin{figure}[ht]
	\centering
	\includegraphics[width=.5\columnwidth]{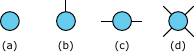}
	\caption{Graphical TN representation of a (a) scalar, (b) vector, (c) matrix and (d) 4-way tensor. \label{fig:TN-vector-matrix-tensor}}
\end{figure}

A \emph{tensor network} (TN) consists of tensors and connections between them that indicate how they should be contracted with each other.
Since the mathematical notation easily gets unwieldy when working with TNs, it is common practice to use graphical representations when discussing such networks.
Figure~\ref{fig:TN-vector-matrix-tensor} shows how scalars, vectors, matrices and tensors can be represented graphically.
A circle or node is used to represent a tensor, and \emph{dangling edges} are used to indicate modes of the tensor.

\begin{figure}[ht]
	\centering
	\includegraphics[width=.5\columnwidth]{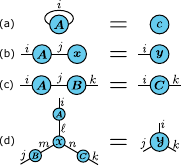}
	\caption{
		TN representation of (a) matrix trace, (b) matrix-vector multiplication, (c) matrix-matrix multiplication and (d) a 3-way Tucker decomposition. 
		\label{fig:TN-simple-contractions}}
\end{figure}

Contraction between two tensors or a tensor with itself can be represented by connecting the appropriate dangling edges.
This is illustrated in Figure~\ref{fig:TN-simple-contractions}.
In mathematical notation, these contractions can be written elementwise as
\begin{enumerate}[(a)]
	\item $\sum_{i} \Abf_{ii} = \trace(\Abf) = c$,
	\item $\sum_{j} \Abf_{ij} \, \xbf_j = \ybf_i$, 
	\item $\sum_{j} \Abf_{ij}\, \Bbf_{jk} = \Cbf_{ik}$,
	\item $\sum_{\ell m n} \Xe_{\ell m n}  \Abf_{i \ell}  \Bbf_{j m}  \Cbf_{k n} = \Ye_{ijk}$. \label{it:Tucker-contraction}
\end{enumerate}

To reduce computational cost when contracting a TN, it is optimal to contract two tensors at a time \citep{pfeifer2014FasterIndentification}.
Determining the optimal contraction order is NP-hard, and developing heuristics for this problem is an active research area \citep{pfeifer2014FasterIndentification, meirom2022OptimizingTensor}.
There are several software packages that compute exact or approximate optimal contraction orders, e.g., NCON \citep{pfeifer2014NCONTensor} for Matlab and the \verb|opt_einsum| package for Python\footnote{Available at \url{https://github.com/dgasmith/opt_einsum}.}.

\begin{figure}[ht]
	\centering
	\includegraphics[width=.4\columnwidth]{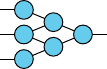}
	\caption{
		Example of a TN matrix.
		\label{fig:TN-matrix}}
\end{figure}

A \emph{TN matrix} is a matrix which is represented by a TN.
The dangling edges in such a network represent either rows or columns of the matrix.
Figure~\ref{fig:TN-matrix} shows an example of a TN matrix with dangling edges pointing to the left representing rows and dangling edges pointing right representing columns.
Suppose there are $N_r$ and $N_c$ dangling edges representing rows and columns, respectively, and let $N = N_r + N_c$.
Let $\Ae \in \Rb^{I_1 \times \cdots \times I_N}$ denote the $N$-way TN representing the matrix, and let $\jbf$ be a length-$N$ vector with the first $N_r$ entries enumerating the modes of $\Ae$ corresponding to rows and the remaining entries enumerating those modes that correspond to columns.
The matrix $\Abf$ represented by $\Ae$ is then given elementwise by
\begin{equation} \label{eq:tensor-network-matrix}
	\Abf(\overline{i_{\jbf(1)} \cdots i_{\jbf(N_r)}}, \overline{i_{\jbf(N_r+1)} \cdots i_{\jbf(N)}}) = \Ae_{i_1 \cdots i_{N}}
\end{equation}
and is of size $\prod_{n=1}^{N_r} I_{\jbf(n)} \times \prod_{m=N_r+1}^{N} I_{\jbf(m)}$.
For example, if the mode dimensions corresponding to the dangling edges in Figure~\ref{fig:TN-matrix} is 100, then the matrix represented by the TN is of size $10^6 \times 100$.

Graphical TN notation does not specify the order in which the modes corresponding to dangling edges should appear when the tensor is represented as a multidimensional array.
A similar ambiguity exists for the TN matrices where it is not clear how the row and column modes should be permuted.
In both cases, any ordering can be used as long the choice is consistent in all computations.

\subsection{Leverage Score Sampling}

In the context of least squares, leverage scores indicate how important the rows of the design matrix are.
By sampling according to this importance metric, it is possible to compute a good approximation (with high probability) to the least squares solution using a random subset of the equations of the full least squares problem.
We provide the necessary preliminaries on leverage score sampling in this section and refer the reader to \citep{mahoney2011RandomizedAlgorithms, woodruff2014SketchingTool} for further details.

\begin{definition} \label{def:leverage-score}
	The $i$th \emph{leverage score} of a matrix $\Abf \in \Rb^{I \times R}$ is defined as $\ell_i(\Abf) \defeq \ebf^{(i)\top} \Abf \Abf^+ \ebf^{(i)}$ for $i \in [I]$, 
	where $\ebf^{(i)}$ is the $i$th canonical basis vector of length $I$.
\end{definition}

\begin{definition} \label{def:leverage-score-sampling-matrix}
	Let $\Abf \in \Rb^{I \times R}$ and let $\pbf \in \Rb^{I}$ be a vector with entries $\pbf(i) = \ell_i(\Abf) / \rank(\Abf)$.
	Then $\pbf$ is the \emph{leverage score distribution} on row indices of $\Abf$.
	Let $f : [J] \rightarrow [I]$ be a random map such that each $f(j)$ is independent and distributed according to $\pbf$.
	The random matrix $\Sbf \in \Rb^{J \times I}$ defined elementwise via
	\begin{equation} \label{eq:sketching-operator-S}
		\Sbf(j, i) \defeq \Ind\{ f(j)=i \} / \sqrt{J \pbf(f(j))}
	\end{equation}
	is the \emph{leverage score sampling matrix} for $\Abf$.
	In \eqref{eq:sketching-operator-S}, $\Ind\{A\}$ is the indicator function which is 1 if the random event $A$ occurs and 0 otherwise.
\end{definition}

The following result is well-known and variants have appeared throughout the literature \citep{drineas2006SamplingAlgorithms, drineas2008RelativeerrorCUR, drineas2011FasterLeast, larsen2022PracticalLeverageBased}.
\begin{theorem} \label{thm:leverage-score-LS}
	Let $\Abf \in \Rb^{I \times R}$ be a matrix and suppose $\Sbf \in \Rb^{J \times I}$ is the leverage score sampling matrix for $\Abf$.
	Moreover, let $\varepsilon, \delta \in (0,1)$, and define $\OPT \defeq \min_{\Xbf} \| \Abf \Xbf - \Ybf \|_\F$ and 
	\begin{equation} \label{eq:X-tilde}
		\tilde{\Xbf} \defeq \argmin_{\Xbf} \| \Sbf \Abf \Xbf - \Sbf \Ybf\|_\F.
	\end{equation}
	If $J \gtrsim R \max(\log(R/\delta), 1/(\varepsilon \delta))$, then with probability at least $1-\delta$ it holds that $\| \Abf \tilde{\Xbf} - \Ybf \|_\F \leq (1+\varepsilon) \OPT$.
\end{theorem}

\section{Proposed Sampling Method} \label{sec:proposed-method}

In this section, we propose an efficient method for sampling rows of a tall-and-skinny TN matrix $\Abf \in \Rb^{I \times R}$ according to its leverage scores.
The method is useful when $I$ is so large that costs and storage on the order $\Omega(I)$ are too expensive, and $R$ is small enough that costs and storage on the order $O(R^3)$ are affordable.
In this setting it is infeasible to sample the rows of $\Abf$ directly via the probability distribution in Definition~\ref{def:leverage-score-sampling-matrix} since
doing so would require computing the pseudoinverse of $\Abf$ which costs $O(I R^2)$.
Moreover, once that pseudoinverse is computed, we would be required to compute, store and sample according to a length-$I$ vector of probabilities.

In view of Definition~\ref{def:leverage-score}, note that
\begin{equation} \label{eq:lev-score-formula}
	\ell_i(\Abf) = \ebf^{(i)\top} \Abf \Abf^+ \ebf^{(i)} = \ebf^{(i)\top} \Abf \Phibf \Abf^\top \ebf^{(i)},
\end{equation}
where $\Phibf \defeq (\Abf^\top \Abf)^+$.
This formulation is the basis for our sampling method.
As we show in Section~\ref{sec:computing-phi}, the matrix $\Phibf$ can be computed efficiently when $\Abf$ is a TN matrix.
Once $\Phibf$ is computed, we can again leverage the TN structure of $\Abf$ to draw samples from its leverage score distribution via the formula in \eqref{eq:lev-score-formula} without forming a length-$I$ probability vector.
We describe how this is done in Section~\ref{sec:drawing-samples-efficiently}.

\subsection{Computing \texorpdfstring{$\Phibf$}{Phi}} \label{sec:computing-phi}

\begin{figure}[ht]
	\centering
	\includegraphics[width=.6\columnwidth]{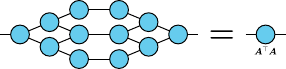}
	\caption{
		Example of how the Gram matrix $\Abf^\top \Abf$ can be computed as a contraction of its TN representation.
		\label{fig:TN-gram}}
\end{figure}

When $\Abf$ is a TN matrix, it is straightforward to represent the Gram matrix $\Abf^\top \Abf$ in such a format as well.
This can be done by linking up the dangling column edges of $\Abf^\top$ with the corresponding dangling row edges of $\Abf$.
For example, suppose $\Abf$ is the TN matrix in Figure~\ref{fig:TN-matrix}.
The graphical representation of $\Abf^\top$ is then a horizontal mirror image of Figure~\ref{fig:TN-matrix}.
The resulting Gram matrix is illustrated in Figure~\ref{fig:TN-gram}.

A dense representation of the Gram matrix can be computed by contracting its TN representation.
For a wide range of TN matrices, this will be much more efficient that computing $\Abf^\top \Abf$ na\"{i}vely by first forming $\Abf$ as a dense matrix and then carrying out a dense matrix-matrix multiplication between $\Abf^\top$ and $\Abf$.
Once $\Abf^\top \Abf$ is computed, its pseudoinverse $\Phibf = (\Abf^\top \Abf)^+$ is affordable to compute via standard methods.

\subsection{Drawing Samples Efficiently} \label{sec:drawing-samples-efficiently}

We now describe how to sample rows of $\Abf$ according to its leverage score distribution without needing to form $\Abf$.
In order to keep the notation simple, we will assume without loss of generality that the modes of the TN $\Ae$ describing the matrix $\Abf$ are arranged so that
\begin{equation*}
	\Abf(\overline{i_1 \cdots i_{N_r}}, \overline{i_{N_r+1} \cdots i_N}) = \Ae_{i_1 \cdots i_N}.
\end{equation*}
Sampling a row of $\Abf$ is therefore equivalent to sampling a multi-index $\overline{i_1 \cdots i_{N_r}} \in [\prod_{n=1}^{N_r} I_n]$ with each $i_n \in [I_n]$.

In a nutshell, our sampling scheme sequentially samples each index $i_1, \ldots, i_{N_r}$ one after another by conditioning on the previously drawn indices.
This allows us to sample from the leverage score distribution while at the same time avoid forming the full probability vector of length $\prod_{n=1}^{N_r} I_n$.

\begin{figure}[ht]
	\centering
	\includegraphics[width=.7\columnwidth]{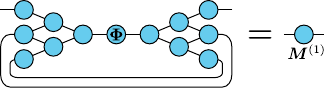}
	\caption{
		Example of TN matrix for computing the distribution of the first index $i_1$.
		\label{fig:TN-sampling-1}}
\end{figure}

In the following we use $\Pb(i_1)$ to denote the probability that the first index is $i_1$.
We use $\Pb(i_n \mid i_1, \ldots, i_{n-1})$ to denote the conditional probability that the $n$th index is $i_n$ given that the previous $n-1$ indices are $i_1, \ldots, i_{n-1}$.
Let $\rho \defeq \rank(\Abf)$. 
From Definition~\ref{def:leverage-score-sampling-matrix} and the formula in \eqref{eq:lev-score-formula} we would like to sample row $\overline{i_1 \cdots i_{N_r}}$ with probability
\begin{equation*}
	\frac{1}{\rho} \, \Abf(\overline{i_1 \cdots i_{N_r}}, :) \, \Phibf \, \Abf(\overline{i_1 \cdots i_{N_r}}, :)^\top.
\end{equation*}
Therefore, we want the first index to be $i_1$ with probability 
\begin{equation} \label{eq:probability-i1}
	\Pb(i_1) = \frac{1}{\rho} \sum_{i_2 \cdots i_{N_r}} \Abf(\overline{i_1 \cdots i_{N_r}}, :) \, \Phibf \, \Abf(\overline{i_1 \cdots i_{N_r}}, :)^\top.
\end{equation}
The matrix $\Abf \Phibf \Abf^\top$ can easily be formulated as a TN matrix by linking up the TNs for $\Abf$ and $\Abf^\top$ with $\Phibf$.
Moreover, the summation in \eqref{eq:probability-i1} in amounts to adding an additional $N_r-1$ contractions to the TN representation of $\Abf \Phibf \Abf^\top$.
All in all, this results in a TN matrix $\Mbf^{(1)} \in \Rb^{I_1 \times I_1}$ with only two dangling edges.
The desired probabilities then lie along the diagonal of this matrix: $\Pb(i_1) = \Mbf^{(1)}_{i_1 i_1}$.
The dense representation of $\Mbf^{(1)}$ can be computed efficiently by contracting the underlying network (it is sufficient to only compute diagonal entries which further reduces the cost).
An example of what this looks like when $\Abf$ is the TN matrix in Figure~\ref{fig:TN-matrix} is illustrated in Figure~\ref{fig:TN-sampling-1}.
The first index $i_1$ is then drawn according to the probability distribution $(\Pb(i_1))_{i_1=1}^{I_1}$.

\begin{figure}[ht]
	\centering
	\includegraphics[width=.7\columnwidth]{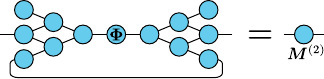}
	\caption{
		Example of TN matrix for computing the distribution of the second index $i_2$ conditionally on the first index $i_1$.
		\label{fig:TN-sampling-2}}
\end{figure}

The next step is to draw subsequent indices conditionally on the previously drawn indices.
Suppose we have drawn $i_1, \ldots, i_{n-1}$.
When sampling the $n$th index we compute the probabilities
\begin{equation} \label{eq:conditional-prob}
	\Pb(i_n \mid i_1, \ldots, i_{n-1}) = \frac{\Pb(i_1, \ldots, i_{n})}{\Pb(i_1, \ldots, i_{n-1})}
\end{equation}
for each $i_n \in [I_n]$.
Since $\Pb(i_1, \ldots, i_{n-1}) = \sum_{i_n} \Pb(i_1, \ldots, i_n)$ it is sufficient to compute the numerator in \eqref{eq:conditional-prob} for each $i_n \in [I_n]$ and then normalize them so that they add up to 1.
Note that
\begin{equation} \label{eq:probability-in}
	\Pb(i_1, \ldots, i_{n}) = 
	\frac{1}{\rho} \sum_{i_{n+1} \cdots i_{N_r}} \Abf(\overline{i_1 \cdots i_{N_r}}, :) \, \Phibf \, \Abf(\overline{i_1 \cdots i_{N_r}}, :)^\top.
\end{equation}
As earlier, this can be represented as a TN by adding $N_r - n$ additional contractions over the indices $i_{n+1}, \ldots, i_{N_r}$ to the TN representation of $\Abf \Phibf \Abf^\top$.
Moreover, since the indices $i_1, \ldots, i_{n-1}$ are fixed, this effectively eliminates the corresponding modes of the underlying tensor representation by reducing their dimensionality to 1.
This results in a TN matrix $\Mbf^{(n)} \in \Rb^{I_n \times I_n}$ with only two dangling edges, with the desired joint probabilities laying along the diagonal: $\Pb(i_1,\ldots,i_n) = \Mbf^{(n)}_{i_n i_n}$ (again, the cost can be reduced by only computing the diagonal elements of $\Mbf^{(n)}$).
Figure~\ref{fig:TN-sampling-2} illustrates what this looks like when we are computing the probability $\Pb(i_1, i_2)$ for the example from Figure~\ref{fig:TN-sampling-1} with $i_1$ fixed.

This sampling procedure is carried out until all indices $i_1,\ldots,i_{N_r}$ have been drawn, which corresponds to drawing a single row index $\overline{i_1 \cdots i_{N_r}}$.
The procedure is then repeated until the desired number of samples has been drawn.
Note that the distribution $(\Pb(i_1))_{i_1=1}^{I_1}$ remains the same for all samples, so in order to speed up the overall sampling procedure it is a good idea to draw $i_1$ for all samples immediately.

The efficiency of the proposed sampling scheme relies on the possibility to efficiently contract the TNs corresponding to the Gram matrix computation $\Abf^\top \Abf$ and the probability distributions in \eqref{eq:probability-i1} and \eqref{eq:conditional-prob}.
For a large class of tensor networks, this will be possible.
We discuss the computational complexity of these contractions for the TN matrices that arise for CP and tensor ring decompositions in Section~\ref{sec:detailed-complexity}.

\section{Application to Tensor Decomposition} \label{sec:application-to-tensor-decomposition}

TNs can be used to express a very wide range of tensor decompositions, including the CP \citep{hitchcock1927ExpressionTensor}, Tucker \citep{tucker1963ImplicationsFactor}, tensor train \citep{oseledets2011TensortrainDecomposition} and tensor ring decompositions \citep{zhao2016TensorRing}.

\begin{figure}[ht]
	\centering
	\includegraphics[width=.7\columnwidth]{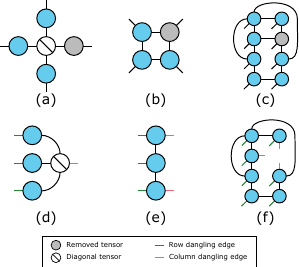}
	\caption{
		The top plots show examples of tensor decompositions: 
		(a) CP decomposition, (b) tensor ring decomposition, and (c) a TN decomposition discovered in \citep{li2020EvolutionaryTopology} via an evolutionary search procedure.
		The bottom plots (d)--(f) show examples of TN matrices that arise when updating the tensor marked in gray in the corresponding top plot.
		A tensor $\Ae$ is said to be \emph{diagonal} if only the elements $\Ae_{i i \cdots i}$ are nonzero.
		\label{fig:TN-tensor-decompositions}}
\end{figure}

Suppose $\Xe \in \Rb^{I_1 \times \cdots I_N}$ is a data tensor that we want to decompose into some TN format consisting of tensors $\Ae^{(1)}, \ldots, \Ae^{(M)}$.
We use $\TN(\Ae^{(1)}, \ldots, \Ae^{(M)}) \in \Rb^{I_1 \times \cdots \times I_N}$ to denote this TN.
This decomposition problem can be formulated as the optimization problem
\begin{equation} \label{eq:tn-optimization-problem}
	\argmin_{\Ae^{(1)}, \ldots, \Ae^{(M)}} \| \Xe - \TN(\Ae^{(1)}, \ldots, \Ae^{(M)}) \|_\F,
\end{equation}
where $\|\cdot\|_\F$ is a straightforward generalization of the matrix Frobenius norm to tensors (see Section~\ref{sec:additional-notation}).
In general, this is a difficult non-convex optimization problem which is hard to solve.
A popular approach to finding an approximate solution to \eqref{eq:tn-optimization-problem} is to iteratively update one tensor $\Ae^{(m)}$ at a time via \emph{alternating least squares} (ALS).
Minimizing the objective in \eqref{eq:tn-optimization-problem} with respect to a single tensor at a time turns it into a \emph{linear} least squares problem:
\begin{equation} \label{eq:tn-optimization-als}
	\argmin_{\Abf^{(m)}} \| \Xbf^{(m)} - \Abf^{\neq m} \Abf^{(m)}\|_\F,
\end{equation}
where $\Xbf^{(m)}$ and $\Abf^{(m)}$ are appropriate matricizations of $\Xe$ and $\Ae^{(m)}$, respectively, and $\Abf^{\neq m}$ is a TN matrix which depends on the tensors $\Ae^{(1)}, \ldots, \Ae^{(m-1)}, \Ae^{(m+1)}, \ldots, \Ae^{(M)}$.
For example, for the CP and Tucker decompositions, $\Abf^{\neq m}$ takes the form of a Khatri--Rao and Kronecker product matrix, respectively \citep{kolda2009TensorDecompositions}.
Figure~\ref{fig:TN-tensor-decompositions} shows examples of three tensor decompositions (top plots) and examples of what a TN matrix $\Abf^{\neq m}$ corresponding to each decomposition looks like (bottom plots).
Algorithm~\ref{alg:TN-ALS} provides a high-level ALS algorithm for computing any TN decomposition.

\begin{algorithm}
	\caption{ALS for TN decomposition}
	\label{alg:TN-ALS}
	\DontPrintSemicolon 
	\KwIn{$\Xe \in \Rb^{I_1 \times \cdots \times I_N}$}
	\KwOut{Factorization tensors $\Ae^{(1)}, \ldots, \Ae^{(M)}$}
	Initialize tensors $\Ae^{(2)}, \ldots, \Abf^{(M)}$ \label{line:tn-als:initialization}\;
	
	\While{termination criteria not met}{
		\For{$m = 1,\ldots, M$}{
			$\Abf^{(m)} = \argmin_{\Abf} \| \Xbf^{(m)} - \Abf^{\neq m} \Abf \|_\F$ \label{line:tn-als:ls} \;
		}	
	}
	\Return{$\Ae^{(1)}, \ldots, \Ae^{(M)}$} \;
\end{algorithm}

ALS is the ``workhorse approach'' to tensor decomposition \citep{kolda2009TensorDecompositions}.
ALS-based methods for various decompositions are implemented in many leading tensor software packages, including Tensor Toolbox \citep{bader2021MATLABTensor} and Tensorlab \citep{vervliet2016Tensorlab} for Matlab and TensorLy \citep{kossaifi2019TensorLyTensor} for Python.
ALS typically works well, especially when the data tensor $\Xe$ is not too large.
However, as previously discussed, each iteration of ALS depends on all entries of $\Xe$.
If $\Xe$ has $N$ modes each of dimension $I$, then $\Xe$ contains $I^N$ entries.
Algorithms for tensor decomposition therefore inherit this exponential dependence on $N$.
In the case of ALS, \emph{each iteration} computed via \eqref{eq:tn-optimization-als} requires access to \emph{all} entries of $\Xe$.
Moreover, for every iteration, it requires forming the matrix $\Abf^{\neq m}$ and then solving a linear system involving this matrix.

As discussed in Section~\ref{sec:related-work}, several previous works have sought to address the curse of dimensionality in ALS algorithms by sampling the least squares problem in \eqref{eq:tn-optimization-als}.
These works develop methods for specific tensor decompositions.
By contrast, our sampling framework described in Section~\ref{sec:proposed-method} can be used to develop a sampling-based decomposition scheme for \emph{any} TN decomposition.
This is done by independently sampling each least squares problem on line~\ref{line:tn-als:ls} in Algorithm~\ref{alg:TN-ALS} according to the leverage scores of the TN matrix $\Abf^{\neq m}$.
The sampling is done by first determining the indices to sample by applying our techniques in Section~\ref{sec:proposed-method}.
Letting $\Sbf$ represent the sampling operation, the next step is to compute $\Sbf \Xbf^{(m)}$ and $\Sbf \Abf^{\neq m}$.
The final step is to solve the sampled least squares problem
\begin{equation} \label{eq:sketched-LS-TN-decomposition}
	\Abf^{(m)} = \argmin_{\Abf} \| \Sbf \Xbf^{(m)} - \Sbf \Abf^{\neq m} \Abf\|_\F.
\end{equation}

In order to make the discussion in Sections~\ref{sec:proposed-method} and \ref{sec:application-to-tensor-decomposition} a bit more concrete, we now apply these ideas to the CP decomposition.

\begin{Example} \label{ex:CP-decomposition}
	Suppose $\Xe \in \Rb^{I_1 \times \cdots \times I_N}$ is an $N$-way tensor and that we want to compute a rank-$R$ CP decomposition of it.
	Figure~\ref{fig:TN-tensor-decompositions}~(a) shows what the CP decomposition looks like when $N=4$.
	Mathematically, for an arbitrary $N$ it takes the form
	\begin{equation} \label{eq:cp-decomposition}
		\Xe_{i_1 \cdots i_N} \approx \sum_{r=1}^R \lambda_r \Abf^{(1)}_{i_1 r} \cdots \Abf^{(N)}_{i_N r},
	\end{equation} 
	where $\lambda_1, \ldots, \lambda_R$ are the elements in the diagonal tensor in the middle of Figure~\ref{fig:TN-tensor-decompositions}~(a), and the other tensors now have two modes and are therefore denoted by $\Abf^{(m)} \in \Rb^{I_m \times R}$ and referred to as \emph{factor matrices}.
	The values of the diagonal elements can be incorporated into the factor matrices, and we therefore assume that they are all equal to 1 without loss of generality.
	The design matrices now take the form 
	\begin{equation*}
		\Abf^{\neq m} = \Abf^{(N)} \odot \cdots \Abf^{(m+1)} \odot \Abf^{(m-1)} \odot \cdots \odot \Abf^{(1)}.
	\end{equation*}
	The corresponding \emph{unfoldings} $\Xbf^{(m)}$ are defined elementwise via
	\begin{equation*}
		\Xbf^{(m)}(\overline{i_1 \cdots i_{m-1} i_{m+1} \cdots i_N}, i_m) = \Xe_{i_1 \cdots i_N}.
	\end{equation*}
	The TN corresponding to the Gram matrix $\Abf^{\neq m \top} \Abf^{\neq m}$ can be efficiently contracted via the well-known formula
	\begin{equation} \label{eq:CP-Gram}
		\Abf^{\neq m \top} \Abf^{\neq m} = \startimes_{j \neq m} \Abf^{(j)\top} \Abf^{(j)},
	\end{equation}
	where $\circledast$ denotes elementwise product.
	Drawing a row index $i \in [\prod_{j \neq m} I_j]$ for $\Abf^{\neq m}$ corresponds to drawing a multi-index $\overline{i_1 \cdots i_{m-1} i_{m+1} \cdots i_N}$ with each $i_j \in [I_j]$.
	It is straightforward to show that the appropriate sampling probabilities in \eqref{eq:probability-i1} and \eqref{eq:probability-in} are now given by
	\begin{equation*}
		\Pb((i_j)_{j \leq n, j \neq m}) = \frac{1}{R} \sum_{r k} \Phibf_{rk} 
		\cdot \Big( \prod_{\substack{j \leq n \\ j \neq m}} \Abf^{(j)}_{i_j r} \Abf^{(j)}_{i_j k} \Big) \Big( \prod_{\substack{j > n \\ j \neq m}} (\Abf^{(j)\top} \Abf^{(j)})_{rk} \Big).
	\end{equation*}
	This formula appears in Lemma~10 of \citep{malik2022MoreEfficient}, but in that formula $\Phibf$ is an approximation of $(\Abf^{\neq m \top} \Abf^{\neq m})^+$.
	This also explains why a different normalization constant is used in \citep{malik2022MoreEfficient}.
	The proof, however, remains the same.
	
\end{Example}

\subsection{Computational Complexity Examples} \label{sec:tn-decomp-complexity}

In this section we give the computational complexity of our proposed approach when used for CP and tensor ring decomposition.
This allows us to compare our work with \citep{malik2022MoreEfficient} which also considers these two decompositions.
Tables~\ref{tab:cp-complexity-comparison} and \ref{tab:tr-complexity-comparison} present the computational complexity of various methods for CP and tensor ring decomposition, respectively, when doing a rank-$R$ decomposition of an $N$-way cubic tensor of size $\Xe \in \Rb^{I \times \cdots \times I}$. 
For a precise definition of the rank $R$ for each decomposition, see \eqref{eq:cp-decomposition} and \eqref{eq:tr-decomposition}.
Our approaches for CP and tensor ring decomposition are called TNS-CP and TNS-TR, respectively, where TNS stands for Tensor Network Sampling. 
The complexities for our methods are computed in Section~\ref{sec:detailed-complexity}.
The complexities for the competing methods are taken from Tables~1 and 2 in \citep{malik2022MoreEfficient}.
For the randomized methods, the complexities reflect the parameter choices required for theoretical performance guarantees; see the respective cited works for further details.

As mentioned in Section~\ref{sec:related-work}, the primary contribution of \citep{malik2022MoreEfficient} was to reduce the dependence on $N$ in the per-iteration cost from exponential to polynomial in ALS for CP and tensor ring decomposition.
Our methods further improve on that cost by reducing the dependence on $R$.
The improvement is due to the avoidance of the complicated recursive sketching step used in \citep{malik2022MoreEfficient}.
Note that all methods other than ours and those by \citep{malik2022MoreEfficient} have an exponential dependence on $N$.

\begin{table}[htb] 
	\centering
	\caption{
		Leading order computational cost (ignoring $\log$ factors) for various CP decomposition methods.
		$\noiter$ is the number of ALS iterations.
		$R$ denotes the CP rank.
		Note that $\tilde{\varepsilon}$ is used as the accuracy parameter instead of $\varepsilon$ for SPALS.
		The reason is that SPALS has additive-error guarantees which are weaker than the relative-error guarantees of the other randomized methods; 
		see Section~\ref{sec:remarks-on-tables} for a discussion.
		\label{tab:cp-complexity-comparison}
	}
	\begin{tabular}{ll}  
		\toprule
		Method    											& Complexity 												\\
		\midrule
		CP-ALS \citep{kolda2009TensorDecompositions}		& $\noiter \cdot N (N+I) I^{N-1} R$ 						\\
		SPALS \citep{cheng2016SPALSFast} 					& $I^N + \noiter \cdot N (N+I) R^{N+1} / \tilde{\varepsilon}^2$ 			\\	
		CP-ARLS-LEV	\citep{larsen2022PracticalLeverageBased}& $\noiter \cdot N ( R + I ) R^N / (\varepsilon \delta)$	\\
		CP-ALS-ES \citep{malik2022MoreEfficient}   			& $\noiter \cdot N^2 R^3 (R + N I / \varepsilon) / \delta$	\\
		\textbf{TNS-CP (our)}								& $\noiter \cdot N^3 I R^3 / (\varepsilon \delta)$ 			\\
		\bottomrule
	\end{tabular}
\end{table}

\begin{table}[htb] 
	\centering
	\caption{
		Leading order computational cost (ignoring $\log$ factors) for various tensor ring decomposition methods. 
		$\noiter$ is the number of ALS iterations.
		The tensor ring rank is $(R,\ldots,R)$.
		\label{tab:tr-complexity-comparison}
	}
	\begin{tabular}{ll}  
		\toprule
		Method    						& Complexity \\
		\midrule
		TR-ALS \citep{zhao2016TensorRing}						& $\noiter \cdot N I^N R^2$ \\
		rTR-ALS \citep{yuan2019RandomizedTensor}							& $N I^N K + \noiter \cdot N K^N R^2$ \\
		TR-SVD \citep{zhao2016TensorRing}							& $I^{N+1} + I^N R^3$ \\
		TR-SVD-Rand \citep{ahmadi-asl2020RandomizedAlgorithms}						& $I^N R^2$ \\
		TR-ALS-Sampled \citep{malik2021SamplingBasedMethod}					& $\noiter \cdot N I R^{2N+2} / (\varepsilon \delta)$	\\
		TR-ALS-ES \citep{malik2022MoreEfficient}	   					& $\noiter \cdot N^3 R^8 ( R + I / \varepsilon)/\delta$	\\
		\textbf{TNS-TR (our)}			& $\noiter \cdot N^3 I R^8 / (\varepsilon \delta)$ \\
		\bottomrule
	\end{tabular}
\end{table}

\section{Experiments} \label{sec:experiments}

We run the same feature extraction and classification experiment as in \citep{malik2022MoreEfficient}.
The experiment considers 7200 images from the COIL-100 dataset \citep{nene1996ColumbiaObject}.
These images depict 100 different objects (e.g., a toy car, an onion, a tomato), each photographed from 72 different angles. 
Each image is 128 by 128 pixels and has three color channels.
We arrange these images into a tensor $\Xe$ of size $7200 \times 128 \times 128 \times 3$ and then decompose it using either a rank-25 CP decomposition or a rank-$(5,5,5,5)$ tensor ring decomposition.
For the CP decomposition, we treat the rows of the factor matrix $\Abf^{(1)} \in \Rb^{7200 \times 25}$ as feature vectors for each image.
We treat 90\% of the images as labeled and use them to classify the images in the remaining 10\% using a $k$-nearest neighbor algorithm with $k=1$ and 10-fold cross validation.
The tensor ring decomposition is used in the same way, but instead of the CP factor matrix we now use the core tensor $\Ae^{(1)} \in \Rb^{5 \times 7200 \times 5}$ corresponding to the first mode and reshape the slices $(\Ae^{(1)}_{:i:})_{i}$ into length-25 feature vectors.
Further details on the experiment setup are given in Section~\ref{sec:further-experiment-details}.

Table~\ref{tab:classification-results} shows the run time, relative decomposition error (Err.) and classification accuracy (Acc.) for a number of different CP and tensor ring decomposition methods.
The algorithms for our TNS-CP and TNS-TR were implemented\footnote{Our code is available at \url{https://github.com/OsmanMalik/TNS}.} by modifying the codes for CP-ALS-ES and TR-ALS-ES by \citep{malik2022MoreEfficient} appropriately.

\begin{table}[ht!]
	\centering
	\caption{
		Run time, decomposition error and classification accuracy in the feature extraction experiment.
		\label{tab:classification-results}
	}
	\begin{tabular}{lrrr}  
		\toprule
		Method    						& Time (s) & Err.\ & Acc.\ (\%) \\
		\midrule
		CP-ALS (Ten.\ Toolbox)			&    48.0 & 0.31 & 99.2 \\
		CPD-ALS (Tensorlab)				&    71.9 & 0.31 & 99.0 \\
		CPD-MINF (Tensorlab)			&   108.6 & 0.33 & 99.7 \\
		CPD-NLS (Tensorlab)				&   112.1 & 0.31 & 92.6 \\
		CP-ARLS-LEV						&    28.6 & 0.32 & 97.7 \\
		CP-ALS-ES						&    29.0 & 0.32 & 98.3 \\
		\textbf{TNS-CP (our)}			&    23.5 & 0.32 & 98.3 \\
		\midrule
		TR-ALS							& 10278.4 & 0.31 & 99.5 \\
		TR-ALS-Sampled					&     5.3 & 0.33 & 98.0 \\
		TR-ALS-ES	   					&    29.7 & 0.33 & 97.2 \\
		\textbf{TNS-TR (our)}			&    14.5 & 0.33 & 97.3 \\
		\bottomrule
	\end{tabular}
\end{table}

All methods achieve roughly the same decomposition error, with the randomized methods typically having a slightly higher error.
All methods except CPD-NLS yield similar classification accuracy. 
It is clear that our methods are faster than CP-ALS-ES and TR-ALS-ES from \citep{malik2022MoreEfficient}, corroborating the improved complexity reported in Tables~\ref{tab:cp-complexity-comparison} and \ref{tab:tr-complexity-comparison}.

\section{Conclusion} \label{sec:conclusion}

We have presented a sampling-based ALS approach for decomposing a tensor into \emph{any} TN format.
The generality of the framework is notable---all previous works on randomized tensor decomposition we are aware of develop methods for specific decompositions.
Additionally, unlike most previous sampling-based ALS algorithms in the literature, our framework makes it possible to do the sampling according to the \emph{exact} leverage scores.
Our approach simplifies and generalizes the scheme developed in \citep{malik2022MoreEfficient}.
Both complexity analyses and numerical experiments confirm the improved run time we are able to achieve compared to \citep{malik2022MoreEfficient} when we use our framework for CP and tensor ring decomposition.

There are many exciting avenues for future research.
One direction is to adapt and implement our method for use in a high-performance distributed memory setting.
Another interesting direction is designing an exact leverage score sampling scheme for general TN matrices $\Abf \in \Rb^{I \times R}$ that avoid the $\Omega(R^3)$ cost that we incur when computing the pseudoinverse $\Phibf = (\Abf^\top \Abf)^+$.
Such methods are known for certain matrices with Kronecker product structure \citep{malik2022FastAlgorithms}, but for general TN matrices such results are not known.

\section*{Acknowledgments}

OAM was supported by the LDRD Program of Lawrence Berkeley National Laboratory under U.S.\ DOE Contract No.\ DE-AC02-05CH11231.
VB was supported in part by the U.S.\ DOE, Office of Science, Office of Advanced Scientific Computing Research, DOE Computational Science Graduate Fellowship under Award No.\ DE-SC0022158.
RM was supported by NSF Grant No.\ 2004235.

\bibliography{library}

\begin{thebibliography}{31}
\providecommand{\natexlab}[1]{#1}
\providecommand{\url}[1]{\texttt{#1}}
\expandafter\ifx\csname urlstyle\endcsname\relax
  \providecommand{\doi}[1]{doi: #1}\else
  \providecommand{\doi}{doi: \begingroup \urlstyle{rm}\Url}\fi

\bibitem[Ahmadi-Asl et~al.(2020)Ahmadi-Asl, Cichocki, Phan, Asante-Mensah,
  Mousavi, Oseledets, and Tanaka]{ahmadi-asl2020RandomizedAlgorithms}
S.~Ahmadi-Asl, A.~Cichocki, A.~H. Phan, M.~G. Asante-Mensah, F.~Mousavi,
  I.~Oseledets, and T.~Tanaka.
\newblock Randomized algorithms for fast computation of low-rank tensor ring
  model.
\newblock \emph{Machine Learning: Science and Technology}, 2020.

\bibitem[Bader et~al.(2021)Bader, Kolda, and {others}]{bader2021MATLABTensor}
B.~W. Bader, T.~G. Kolda, and {others}.
\newblock {{MATLAB Tensor Toolbox}}, {{Version}} 3.2.1.
\newblock Available online at https://www.tensortoolbox.org, 2021.

\bibitem[Cheng et~al.(2016)Cheng, Peng, Liu, and Perros]{cheng2016SPALSFast}
D.~Cheng, R.~Peng, Y.~Liu, and I.~Perros.
\newblock {SPALS}: Fast alternating least squares via implicit leverage scores
  sampling.
\newblock In \emph{Advances In Neural Information Processing Systems}, pages
  721--729, 2016.

\bibitem[Cichocki et~al.(2016)Cichocki, Lee, Oseledets, Phan, Zhao, and
  Mandic]{cichocki2016TensorNetworksPart1}
A.~Cichocki, N.~Lee, I.~Oseledets, A.-H. Phan, Q.~Zhao, and D.~P. Mandic.
\newblock Tensor networks for dimensionality reduction and large-scale
  optimization: Part 1 low-rank tensor decompositions.
\newblock \emph{Foundations and Trends in Machine Learning}, 9\penalty0
  (4-5):\penalty0 249--429, 2016.

\bibitem[Cichocki et~al.(2017)Cichocki, Phan, Zhao, Lee, Oseledets, Sugiyama,
  and Mandic]{cichocki2017TensorNetworksPart2}
A.~Cichocki, A.-H. Phan, Q.~Zhao, N.~Lee, I.~Oseledets, M.~Sugiyama, and D.~P.
  Mandic.
\newblock Tensor networks for dimensionality reduction and large-scale
  optimization: Part 2 applications and future perspectives.
\newblock \emph{Foundations and Trends in Machine Learning}, 9\penalty0
  (6):\penalty0 431--673, 2017.

\bibitem[Drineas et~al.(2006)Drineas, Mahoney, and
  Muthukrishnan]{drineas2006SamplingAlgorithms}
P.~Drineas, M.~W. Mahoney, and S.~Muthukrishnan.
\newblock Sampling algorithms for $\ell_2$ regression and applications.
\newblock In \emph{Proceedings of the Seventeenth Annual {ACM}-{SIAM} Symposium
  on Discrete Algorithm}, pages 1127--1136, 2006.

\bibitem[Drineas et~al.(2008)Drineas, Mahoney, and
  Muthukrishnan]{drineas2008RelativeerrorCUR}
P.~Drineas, M.~W. Mahoney, and S.~Muthukrishnan.
\newblock Relative-error {{CUR}} matrix decompositions.
\newblock \emph{{SIAM} Journal on Matrix Analysis and Applications},
  30\penalty0 (2):\penalty0 844--881, 2008.

\bibitem[Drineas et~al.(2011)Drineas, Mahoney, Muthukrishnan, and
  Sarl{\'o}s]{drineas2011FasterLeast}
P.~Drineas, M.~W. Mahoney, S.~Muthukrishnan, and T.~Sarl{\'o}s.
\newblock Faster least squares approximation.
\newblock \emph{Numerische Mathematik}, 117\penalty0 (2):\penalty0 219--249,
  2011.

\bibitem[Fahrbach et~al.(2022)Fahrbach, Fu, and
  Ghadiri]{fahrbach2022SubquadraticKronecker}
M.~Fahrbach, T.~Fu, and M.~Ghadiri.
\newblock Subquadratic {K}ronecker regression with applications to tensor
  decomposition.
\newblock \emph{Preprint arXiv:2209.04876}, 2022.

\bibitem[Hitchcock(1927)]{hitchcock1927ExpressionTensor}
F.~L. Hitchcock.
\newblock The expression of a tensor or a polyadic as a sum of products.
\newblock \emph{Journal of Mathematics and Physics}, 6\penalty0 (1-4):\penalty0
  164--189, 1927.

\bibitem[Ji et~al.(2019)Ji, Wang, Li, and Liu]{ji2019SurveyTensor}
Y.~Ji, Q.~Wang, X.~Li, and J.~Liu.
\newblock A survey on tensor techniques and applications in machine learning.
\newblock \emph{IEEE Access}, 7:\penalty0 162950--162990, 2019.

\bibitem[Kolda and Bader(2009)]{kolda2009TensorDecompositions}
T.~G. Kolda and B.~W. Bader.
\newblock Tensor decompositions and applications.
\newblock \emph{SIAM Review}, 51\penalty0 (3):\penalty0 455--500, Aug. 2009.

\bibitem[Kossaifi et~al.(2019)Kossaifi, Panagakis, Anandkumar, and
  Pantic]{kossaifi2019TensorLyTensor}
J.~Kossaifi, Y.~Panagakis, A.~Anandkumar, and M.~Pantic.
\newblock {TensorLy}: Tensor learning in {P}ython.
\newblock \emph{Journal of Machine Learning Research}, 20\penalty0
  (26):\penalty0 1--6, 2019.

\bibitem[Larsen and Kolda(2022)]{larsen2022PracticalLeverageBased}
B.~W. Larsen and T.~G. Kolda.
\newblock Practical leverage-based sampling for low-rank tensor decomposition.
\newblock \emph{SIAM Journal on Matrix Analysis and Applications}, 43\penalty0
  (3):\penalty0 1488--1517, 2022.

\bibitem[Li and Sun(2020)]{li2020EvolutionaryTopology}
C.~Li and Z.~Sun.
\newblock Evolutionary topology search for tensor network decomposition.
\newblock In \emph{Proceedings of the 37th International Conference on Machine
  Learning}, volume 119, pages 5947--5957. PMLR, 2020.

\bibitem[Ma and Solomonik(2022)]{ma2022CostefficientGaussian}
L.~Ma and E.~Solomonik.
\newblock Cost-efficient {G}aussian tensor network embeddings for
  tensor-structured inputs.
\newblock \emph{Preprint arXiv:2205.13163}, 2022.

\bibitem[Mahoney(2011)]{mahoney2011RandomizedAlgorithms}
M.~W. Mahoney.
\newblock Randomized algorithms for matrices and data.
\newblock \emph{Foundations and Trends in Machine Learning}, 3\penalty0
  (2):\penalty0 123--224, 2011.

\bibitem[Malik(2022)]{malik2022MoreEfficient}
O.~A. Malik.
\newblock More efficient sampling for tensor decomposition with worst-case
  guarantees.
\newblock In \emph{Proceedings of the 39th International Conference on Machine
  Learning}, volume 162, pages 14887--14917. PMLR, 2022.

\bibitem[Malik and Becker(2021)]{malik2021SamplingBasedMethod}
O.~A. Malik and S.~Becker.
\newblock A sampling-based method for tensor ring decomposition.
\newblock In \emph{Proceedings of the 38th International Conference on Machine
  Learning}, volume 139, pages 7400--7411. PMLR, 2021.

\bibitem[Malik et~al.(2022)Malik, Xu, Cheng, Becker, Doostan, and
  Narayan]{malik2022FastAlgorithms}
O.~A. Malik, Y.~Xu, N.~Cheng, S.~Becker, A.~Doostan, and A.~Narayan.
\newblock Fast algorithms for monotone lower subsets of {K}ronecker least
  squares problems.
\newblock \emph{Preprint arXiv:2209.05662}, 2022.

\bibitem[Meirom et~al.(2022)Meirom, Maron, Mannor, and
  Chechik]{meirom2022OptimizingTensor}
E.~Meirom, H.~Maron, S.~Mannor, and G.~Chechik.
\newblock Optimizing tensor network contraction using reinforcement learning.
\newblock In \emph{Proceedings of the 39th International Conference on Machine
  Learning}, volume 162, pages 15278--15292. PMLR, 2022.

\bibitem[Nene et~al.(1996)Nene, Nayar, and Murase]{nene1996ColumbiaObject}
S.~A. Nene, S.~K. Nayar, and H.~Murase.
\newblock Columbia object image library ({COIL}-100).
\newblock Technical Report CUCS-006-96, Columbia University, 1996.

\bibitem[Oseledets(2011)]{oseledets2011TensortrainDecomposition}
I.~V. Oseledets.
\newblock Tensor-train decomposition.
\newblock \emph{SIAM Journal on Scientific Computing}, 33\penalty0
  (5):\penalty0 2295--2317, 2011.

\bibitem[Papalexakis et~al.(2016)Papalexakis, Faloutsos, and
  Sidiropoulos]{papalexakis2016TensorsData}
E.~E. Papalexakis, C.~Faloutsos, and N.~D. Sidiropoulos.
\newblock Tensors for data mining and data fusion: Models, applications, and
  scalable algorithms.
\newblock \emph{ACM Transactions on Intelligent Systems and Technology
  ({TIST})}, 8\penalty0 (2):\penalty0 1--44, 2016.

\bibitem[Pfeifer et~al.(2014{\natexlab{a}})Pfeifer, Evenbly, Singh, and
  Vidal]{pfeifer2014NCONTensor}
R.~N. Pfeifer, G.~Evenbly, S.~Singh, and G.~Vidal.
\newblock {NCON}: A tensor network contractor for {MATLAB}.
\newblock \emph{Preprint arXiv:1402.0939}, 2014{\natexlab{a}}.

\bibitem[Pfeifer et~al.(2014{\natexlab{b}})Pfeifer, Haegeman, and
  Verstraete]{pfeifer2014FasterIndentification}
R.~N.~C. Pfeifer, J.~Haegeman, and F.~Verstraete.
\newblock Faster identification of optimal contraction sequences for tensor
  networks.
\newblock \emph{Phys. Rev. E}, 90:\penalty0 033315, 2014{\natexlab{b}}.

\bibitem[Tucker(1963)]{tucker1963ImplicationsFactor}
L.~R. Tucker.
\newblock Implications of factor analysis of three-way matrices for measurement
  of change.
\newblock \emph{Problems in measuring change}, 15\penalty0 (122-137):\penalty0
  3, 1963.

\bibitem[Vervliet et~al.(2016)Vervliet, Debals, Sorber, Van~Barel, and
  De~Lathauwer]{vervliet2016Tensorlab}
N.~Vervliet, O.~Debals, L.~Sorber, M.~Van~Barel, and L.~De~Lathauwer.
\newblock Tensorlab 3.0.
\newblock Available online at https://www.tensorlab.net, Mar. 2016.

\bibitem[Woodruff(2014)]{woodruff2014SketchingTool}
D.~P. Woodruff.
\newblock Sketching as a tool for numerical linear algebra.
\newblock \emph{Foundations and Trends in Theoretical Computer Science},
  10\penalty0 (1-2):\penalty0 1--157, 2014.

\bibitem[Yuan et~al.(2019)Yuan, Li, Cao, and Zhao]{yuan2019RandomizedTensor}
L.~Yuan, C.~Li, J.~Cao, and Q.~Zhao.
\newblock Randomized tensor ring decomposition and its application to
  large-scale data reconstruction.
\newblock In \emph{{IEEE} International Conference on Acoustics, Speech and
  Signal Processing}, pages 2127--2131. {IEEE}, 2019.

\bibitem[Zhao et~al.(2016)Zhao, Zhou, Xie, Zhang, and
  Cichocki]{zhao2016TensorRing}
Q.~Zhao, G.~Zhou, S.~Xie, L.~Zhang, and A.~Cichocki.
\newblock Tensor ring decomposition.
\newblock \emph{Preprint arXiv:1606.05535}, 2016.

\end{thebibliography}

\clearpage
\appendix

\section{Additional Notation} \label{sec:additional-notation}

Let $\Abf \in \Rb^{I \times J}$ and $\Bbf \in \Rb^{K \times L}$ be two matrices.
Their \emph{Kronecker product} is denoted by $\Abf \otimes \Bbf \in \Rb^{IK \times JL}$ and defined as
\begin{equation*}
	\Abf \otimes \Bbf \defeq 
	\begin{bmatrix}
		\Abf_{11} \Bbf & \Abf_{12} \Bbf & \cdots & \Abf_{1J} \Bbf \\
		\Abf_{21} \Bbf & \Abf_{22} \Bbf & \cdots & \Abf_{2J} \Bbf \\
		\vdots & \vdots & & \vdots \\
		\Abf_{I1} \Bbf & \Abf_{I2} \Bbf & \cdots & \Abf_{IJ} \Bbf \\
	\end{bmatrix}.
\end{equation*}

Suppose now that $\Abf \in \Rb^{I \times J}$ and $\Bbf \in \Rb^{K \times J}$, i.e., $\Abf$ and $\Bbf$ have the same number of columns.
The \emph{Khatri--Rao product} of $\Abf$ and $\Bbf$ is denoted by $\Abf \odot \Bbf \in \Rb^{IK \times J}$ and defined as
\begin{equation*}
	\Abf \odot \Bbf \defeq 
	\begin{bmatrix}
		\Abf_{:1} \otimes \Bbf_{:1} & \Abf_{:2} \otimes \Bbf_{:2} & \cdots & \Abf_{:J} \Bbf_{:J}
	\end{bmatrix}.
\end{equation*}
The Khatri--Rao product is sometimes called the \emph{columnwise Kronecker product} for obvious reasons.

Now, suppose that $\Abf \in \Rb^{I \times J}$ and $\Bbf \in \Rb^{I \times J}$ are of the same size.
The \emph{Hadamard product}, or elementwise product, of $\Abf$ and $\Bbf$ is denoted by $\Abf \circledast \Bbf \in \Rb^{I \times J}$ and defined elementwise in the obvious way:
\begin{equation*}
	(\Abf \circledast \Bbf)_{ij} = \Abf_{ij} \Bbf_{ij}.
\end{equation*}

The \emph{Frobenius norm} is a standard matrix norm that can be easily extended to tensors.
For a tensor $\Xe \in \Rb^{I_1 \times \cdots \times I_N}$ it is defined as
\begin{equation*}
	\|\Xe\|_\F \defeq \sqrt{\sum_{i_1 = 1}^{I_1} \cdots \sum_{i_N = 1}^{I_N} \Xe_{i_1 \cdots i_N}^2}.
\end{equation*}

\section{Detailed Complexity Analysis} \label{sec:detailed-complexity}

Here we provide further details on the complexity analysis in Section~\ref{sec:tn-decomp-complexity}.

\subsection{CP Decomposition}

Suppose we are decomposing an $N$-way cubic tensor of size $\Xe \in \Rb^{I \times \cdots \times I}$ using a rank-$R$ CP decomposition.
As mentioned in Example~\ref{ex:CP-decomposition}, the entries of the diagonal tensor (see Figure~\ref{fig:TN-tensor-decompositions}~(a)) can be treated as ones without loss of generality.
The update in \eqref{eq:sketched-LS-TN-decomposition} therefore involves updating one of the CP factor matrices $\Abf^{(m)}$, $m \in [N]$.

\paragraph{Computing $\Phibf$}
Computing the Gram matrix $\Abf^{\neq m \top} \Abf^{\neq m}$ via \eqref{eq:CP-Gram} costs $O(NIR^2)$.
When doing this computation, we store all the Gram matrices $\Abf^{(j) \top} \Abf^{(j)}$ for no additional computational cost for later use.
Computing the pseudoinverse $\Phibf = (\Abf^{\neq m \top} \Abf^{\neq m})^+$ costs an additional $O(R^3)$.
The total cost of computing $\Phibf$ is therefore $O(NIR^2 + R^3)$.

\paragraph{Drawing $J$ Samples}
The cost is the same as the sampling cost for CP-ALS-ES in \citep{malik2022MoreEfficient}, but without the one-time costs since $\Phibf$ and the Gram matrices $\Abf^{(j)\top} \Abf^{(j)}$, $j \in [N] \setminus \{ m \}$, have already been computed.
The cost is therefore $O(J R^2 N^2 I)$; see \citep[Sec.~C.1]{malik2022MoreEfficient} for details.

\paragraph{Solving Sampled Least Squares Problem}
The cost of forming the sampled matrices $\Sbf \Abf^{\neq m}$ and $\Sbf \Xbf^{\neq m}$, and solving the associated linear system \eqref{eq:sketched-LS-TN-decomposition} via direct methods, is the same as for CP-ALS-ES in \citep{malik2022MoreEfficient}, namely $O(JR (N + R + I))$; see \citep[Sec.~C.1]{malik2022MoreEfficient} for details.

\paragraph{Adding It All Up}
The costs above are for solving \emph{one} least squares problem of the form \eqref{eq:sketched-LS-TN-decomposition}. 
Each ALS iteration requires solving $N$ such systems. 
Consequently, the per-iteration cost for our proposed method is
\begin{equation} \label{eq:cp-complexity-J}
	O(J R^2 N^3 I),
\end{equation}
where we have used the fact that the sample size $J$ must satisfy $J > R$ (see Theorem~\ref{thm:leverage-score-LS}) to simplify the complexity expression.
Since $\Abf^{\neq m}$ has $R$ columns, the bound on $J$ from Theorem~\ref{thm:leverage-score-LS} is $J \gtrsim R \max(\log(R/\delta), 1/(\varepsilon\delta))$.
Inserting this into \eqref{eq:cp-complexity-J} gives the per-iteration cost
\begin{equation*}
	O(R^3 N^3 I \log(R/\delta) / (\varepsilon\delta)) = \tilde{O}(R^3 N^3 I / (\varepsilon\delta))
\end{equation*}
which is what we report in Table~\ref{tab:cp-complexity-comparison}.

\subsection{Tensor Ring Decomposition} \label{sec:detailed-complexity-tr}

Suppose that we are decomposing an $N$-way cubic tensor of size $\Xe \in \Rb^{I \times \cdots \times I}$ using a tensor ring decomposition with all ranks set to $R$.
The decomposition takes the form
\begin{equation} \label{eq:tr-decomposition}
	\Xe_{i_1 \cdots i_N} \approx \sum_{r_1 \cdots r_N} \Ae^{(1)}_{r_N i_1 r_1} \Ae^{(2)}_{r_1 i_2 r_2} \cdots \Ae^{(N)}_{r_{N-1} i_N r_N},
\end{equation}
where each $\Ae^{(m)}$ is now a 3-way tensor of size $R \times I \times R$, and where each summation index $r_n$ goes from 1 to $R$.
See \citep{zhao2016TensorRing} for further details on the tensor ring decomposition.

\begin{figure}[htb]
	\centering
	\includegraphics[width=.7\columnwidth]{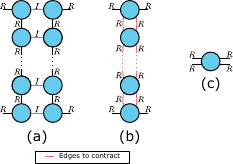}
	\caption{
		Illustration of contraction order of TN for $\Abf^{\neq m \top} \Abf^{\neq m}$ for tensor ring decomposition.
		The variables $I$ and $R$ indicate the dimension of each edge.
		First we eliminate the large dimensions $I$ by contracting over the red edges in (a).
		Then we contract the resulting chain of tensors in sequence as shown in (b).
		The final output is the $R^2 \times R^2$ matrix shown in (c). 
		\label{fig:TR-gram-contraction}}
\end{figure}

\paragraph{Computing $\Phibf$}
For the tensor ring decomposition, the TN representation of $\Abf^{\neq m}$ is shown in Figure~\ref{fig:TN-tensor-decompositions} for the case $N=4$.
Figure~\ref{fig:TR-gram-contraction} illustrates how we contract the TN corresponding to the Gram matrix $\Abf^{\neq m \top} \Abf^{\neq m}$ for an arbitrary $N$.
Each edge we eliminate in subplot (a) costs $O(IR^4)$.
Since there are $N-1$ pairs to contract, this brings the total cost of the step in (a) to $O(NIR^4)$.
Each step in the sequence of contractions in (b) costs $O(R^6)$.
Since $N-2$ such contractions are required, the total cost of step (b) is $O(NR^6)$.
Once the matrix $\Abf^{\neq m \top} \Abf^{\neq m}$ in (c) is computed, it costs $O(R^6)$ to compute its pseudoinverse, i.e., less than the cost of the step in (b).
Adding up the costs for the steps in (a) and (b) and the cost of the pseudoinverse brings the total cost for computing $\Phibf$ to $O(N (R^6 + IR^4))$.

\paragraph{Drawing $J$ Samples}
The cost is the same as the sampling cost for TR-ALS-ES in \citep{malik2022MoreEfficient}, but without the one-time costs in that paper\footnote{%
	The Gram matrices $\Gbf_{[2]}^{(j)\top} \Gbf_{[2]}^{(j)}$ in \citep{malik2022MoreEfficient} are precisely what we compute when we do the contractions in Figure~\ref{fig:TR-gram-contraction}~(a). 
	These can be saved during that step, which is why the one-time costs of computing those Gram matrices is already accounted for from the computation of $\Phibf$.%
}.
The cost is therefore $O(J N^2 I R^6)$; see \citep[Sec.~C.2]{malik2022MoreEfficient} for details.

\paragraph{Solving Sampled Least Squares Problem}
The cost of forming the sampled matrices $\Sbf \Abf^{\neq m}$ and $\Sbf \Xbf^{(m)}$, and solving the associated linear system \eqref{eq:sketched-LS-TN-decomposition} via direct methods, is the same as for TR-ALS-ES in \citep{malik2022MoreEfficient}, namely $O(J (NR^3 + R^4 + IR^2))$.

\paragraph{Adding It All Up}
The costs above are for solving \emph{one} least squares problem of the form \eqref{eq:sketched-LS-TN-decomposition}.
Each ALS iteration requires solving $N$ such systems.
Consequently, the per-iteration cost for our proposed method is
\begin{equation} \label{eq:tr-complexity-J}
	O(J N^3 I R^6).
\end{equation}
Since $\Abf^{\neq m}$ has $R^2$ columns, the bound on $J$ in Theorem~\ref{thm:leverage-score-LS} is $J \gtrsim R^2 \max(\log(R^2/\delta), 1/(\varepsilon\delta))$.
Inserting this into \eqref{eq:tr-complexity-J} gives the per-iteration cost
\begin{equation*}
	O(N^3 I R^8 \log(R^2/\delta) / (\varepsilon \delta)) = \tilde{O}(N^3 I R^8 / (\varepsilon \delta)),
\end{equation*}
which is what we report in Table~\ref{tab:tr-complexity-comparison}.

\subsection{Some Remarks on Tables~\ref{tab:cp-complexity-comparison} and \ref{tab:tr-complexity-comparison}} \label{sec:remarks-on-tables}

All methods in Table~\ref{tab:cp-complexity-comparison} except CP-ALS are randomized. 
The theoretical guarantees for CP-ARLS-LEV, CP-ALS-ES and TNS-CP are all of the same type:
They provide relative-error guarantees of the form in Theorem~\ref{thm:leverage-score-LS} for each least squares solve on line~\ref{line:tn-als:ls} in Algorithm~\ref{alg:TN-ALS} when it is adapted to CP decomposition.
For these three methods, the meaning of the variables $\varepsilon$ and $\delta$ are therefore identical.
SPALS, by contrast, has weaker additive-error guarantees. 
Expressed using the same notation as in Theorem~\ref{thm:leverage-score-LS}, their guarantees take the form 
\begin{equation} \label{eq:additive-error}
	\| \Abf \tilde{\Xbf} - \Ybf \|_\F^2 \leq \OPT^2 + \tilde{\varepsilon} \|\Ybf\|_\F^2,
\end{equation}
where the statement in \eqref{eq:additive-error} holds with probability at least $1-\delta$ when the number of samples is large enough.
The constant $\delta$ has the same meaning as for the other three randomized methods, but $\tilde{\varepsilon}$ has a different meaning than $\varepsilon$: 
The latter is a \emph{relative-error} accuracy parameter while the former is an \emph{additive-error} accuracy parameter.

For the three methods TR-ALS-Sampled, TR-ALS-ES and TNS-TR in Table~\ref{tab:tr-complexity-comparison}, the parameters $\delta$ and $\varepsilon$ are all used in the sense of Theorem~\ref{thm:leverage-score-LS}.

The number of iterations, denoted by $\noiter$, required to reach a certain accuracy or fulfill certain termination criteria may differ between the various methods.
For example, the deterministic ALS methods---by virtue of being exact and non-random---may require fewer iterations to reach a certain accuracy.
Similarly differences may exist between the various randomized ALS methods, and between the different decomposition types.

\section{Further Details on Experiment} \label{sec:further-experiment-details}

\subsection{Competing Methods}

We provide some further details on the methods we compare against below.
\begin{itemize}
	\item \textbf{CP-ALS.} This is an implementation of the standard ALS method for CP decomposition without any randomization.
	It comes with the Tensor Toolbox for Matlab \citep{bader2021MATLABTensor} which can be downloaded at \url{https://www.tensortoolbox.org}.
	
	\item \textbf{CPD-ALS, -MINF, -NLS.} These methods use three different minimization approaches to compute the CP decomposition. 
	They are all available in the Tensorlab package for Matlab \citep{vervliet2016Tensorlab} which can be downloaded at \url{https://www.tensorlab.net/}.
	
	\item \textbf{CP-ARLS-LEV.} This the method proposed in \citep{larsen2022PracticalLeverageBased}.
	We use the implementation by \citet{malik2022MoreEfficient} which is available at \url{https://github.com/OsmanMalik/TD-ALS-ES}.
	
	\item \textbf{CP-ALS-ES.} This is the method for CP decomposition proposed in \citep{malik2022MoreEfficient}.
	We use the code available at \url{https://github.com/OsmanMalik/TD-ALS-ES}.
	
	\item \textbf{TR-ALS.} This is a standard ALS methods for tensor ring decomposition, which was proposed in \citep{zhao2016TensorRing}.
	We use the implementation by \citet{malik2021SamplingBasedMethod} which is available at \url{https://github.com/OsmanMalik/tr-als-sampled}.
	
	\item \textbf{TR-ALS-Sampled.} This is the sampling-based approach for tensor ring decomposition proposed in \citep{malik2021SamplingBasedMethod}.
	We use the code available at \url{https://github.com/OsmanMalik/tr-als-sampled}.
	
	\item \textbf{TR-ALS-ES.} This is the method for tensor ring decomposition proposed in \citep{malik2022MoreEfficient}.
	We use the code available at \url{https://github.com/OsmanMalik/TD-ALS-ES}.
\end{itemize}

\subsection{Parameter Choices}

We sample $J=2000$ rows in all ALS subproblems for the sampling-based CP decomposition methods (CP-ARLS-LEV, CP-ALS-ES and TNS-CP).
For the sampling-based tensor ring methods (TR-ALS-Sampled, TR-ALS-ES and TNS-TR) we use $J=1000$ samples in the ALS subproblems.
Both CP-ALS-ES and TR-ALS-ES require an intermediate embedding via a recursive sketching procedure.
For both methods, we use an embedding dimension of 10000 for these intermediate steps, as suggested in \citep{malik2022MoreEfficient}.

\subsection{Hardware/Software and Dataset}

The experiments are run in Matlab R2021b on a laptop computer with an Intel Core i7-1185G7 CPU and 32 GB of RAM.

The COIL-100 dataset \citep{nene1996ColumbiaObject} is available for download at \url{https://www.cs.columbia.edu/CAVE/software/softlib/coil-100.php}.

\end{document}